\documentclass[12pt]{article}
\usepackage[latin5]{inputenc}
\usepackage{amssymb}
\usepackage{amsmath}

\title{ Multipoint Normal Differential Operators of Second Order}

\author {Erdal UNLUYOL ,
\and Elif OTKUN ÇEVİK
\and and Zameddin I. ISMAILOV  \thanks {Karadeniz Technical University, Art and Science Faculty, Department  of {\rm \; \; \; \; \; \; \; \;} Mathematics, 61080 Trabzon, eunluyol@yahoo.com, zameddin@yahoo.com} }

\date{}

\begin{document}

\maketitle

\begin{abstract} In this work it is described all normal extensions of a multipoint minimal operator generated by linear multipoint differential-operator expresssion for second order in the Hilbert space of vector functions in terms of boundary values at the endpoints of the infinitely many subintervals. Finally, a spectrum structure of such extensions has been investigated.
\end{abstract}

\section{\textbf{Introduction}}

 It is known that traditional infinite direct sum of Hilbert spaces $H_{n} ,{\rm \; }n\ge 1$ and infinite direct sum of operators $A_{n} $ in $H_{n} ,{\rm \; }n\ge 1$ are define as
\[ H=\mathop{\oplus }\limits_{n=1}^{\infty } H_{n} =\left\{u=\left(u_{n} \right):u_{n} \in H_{n} ,{\rm \; }n\ge 1{\rm \; and\; }\sum _{n=1}^{\infty }\left\| u_{n} \right\| _{H_{n} }^{2} <+\infty  \right\}, \]
\begin{multline*}
A=\mathop{\oplus }\limits_{n=1}^{\infty } A_{n}, D(A)=\{u=(u_{n})\in H:u_{n} \in D(A_{n}),{\rm \; }n\ge 1 {\rm \; and }\\ Au=\left(A_{n} u_{n} \right)\in H \}.
\end{multline*}
 A linear space $H$ is a Hilbert space with norm corresponding to inner product
\[\left(u,v\right)_{H} =\sum _{n=1}^{\infty }\left(u_{n} ,v_{n} \right)_{H_{n} }  ,{\rm \; }u,v\in H \cite{Dun}.\]

 The general theory of linear closed operators in Hilbert spaces and its applications to physical problems have been investigated by many mathematicians ( for example, see \cite{Dun}).

 However, many physical problems support to study a theory of linear operators in direct sums in Hilbert spaces ( for example, see \cite{Tim}-\cite{Ism} and references in them )today.

 Notice that a detail analysis of normal subspaces and operators in Hilbert spaces has been studied in \cite{Cod} (see references in it ).

 In this work in first section a connection between multipoint and twopoints normal operators are investigated.

 In second section all normal extensions of multipoint formally normal operators are described in terms of boundary values in the endpoints of the infinitely many subintervals. Furthermore, a spectrum structure has been researched of such extensions.

\section{\textbf{The Minimal and Maximal Operators}}

\indent Along of this work $\left(a_{n} \right)$ and $\left(b_{n} \right)$ will be sequences of real numbers such that
\[-\infty <a_{n} <b_{n} \le a_{n+1} <\cdots <+\infty ,\]
$H_{n} $ is any Hilbert space, $\Delta _{n} =\left(a_{n} ,b_{n} \right),{\rm \; }L_{n}^{2} =L^{2} \left(H_{n} ,\Delta _{n} \right),{\rm \; }L^{2} =\mathop{\oplus }\limits_{n=1}^{\infty } L^{2} \left(H_{n} ,\Delta _{n} \right),\\{\rm \; }\left(\cdot ,\cdot \right)_{H_{n} } =\left(\cdot ,\cdot \right)_{n} $, $n\ge 1$, $W_{2}^{2} =\mathop{\oplus }\limits_{n=1}^{\infty } W_{2}^{2} \left(H_{n} ,\Delta _{n} \right)$, $\mathop{W_{2}^{2} }\limits^{0} =\mathop{\oplus }\limits_{n=1}^{\infty } \mathop{W_{2}^{2} }\limits^{0} \left(H_{n} ,\Delta _{n} \right),{\rm \; }H=\mathop{\oplus }\limits_{n=1}^{\infty } H_{n} ,{\rm \; }cl\left(T\right)$-closure of the operator $T,{\rm \; }E$ is an identity operator in corresponding spaces. $l\left(\cdot \right)$ is a linear multipoint differential-operator expression for second order in $L^{2} $ in the following form

\begin{equation} \label{GrindEQ__2_1_}
l\left(u\right)=\left(l_{n} \left(u_{n} \right)\right)
\end{equation}
and for each $n\ge 1$
\begin{equation} \label{GrindEQ__2_2_}
l_{n} \left(u_{n} \right)=-u''_{n} +iA_{n} u_{n} ,
\end{equation}
where $A_{n} :D\left(A_{n} \right)\subset H_{n} \to H_{n} $ is a linear positive defined selfadjoint operator in $H_{n} $.

 It is clear that formally adjoint expression to \eqref{GrindEQ__2_2_} in the Hilbert space $L_{n}^{2} $ is in the form

\begin{equation} \label{GrindEQ__2_3_}
l_{n}^{+} \left(v_{n} \right)=-v''_{n} -iA_{n}^{*} v_{n} ,{\rm \; }n\ge 1.
\end{equation}

 Define an operator $L'_{n0} $ on the dense manifold of vector functions $D'_{n0} $ in $L_{n}^{2} $,
\begin{multline*}
D'_{n0} :=\{{u_{n} \in L_{n}^{2}} :u_{n} =\sum _{k=1}^{m}\phi _{k} f_{k} ,{\rm \; }\phi _{k} \in C_{0}^{\infty } ({\Delta _{n}}),{\rm \;} f_{k} \in D(A_{n}),\\
k=1,2,\cdots ,m;{\rm \; }m\in {\mathbb N} \}
\end{multline*}
\noindent as $L'_{n0} u_{n} :=l_{n} \left(u_{n} \right),{\rm \; }n\ge 1$.

 Since the operator $A_{n} >0$, then from the relation

\[Im\left(L'_{n0} u_{n} ,u_{n} \right)_{L_{n}^{2} } =\left(A_{n} u_{n} ,u_{n} \right)_{L_{n}^{2} } \ge 0,{\rm \; }u_{n} \in D'_{n0} ,{\rm \; }n\ge 1\]
implies that $L'_{n0} $ is a dissipative in $L_{n}^{2} ,{\rm \; }n\ge 1$. Hence the operator $L'_{n0} $ has a closure in $L_{n}^{2} ,{\rm \; }n\ge 1$. The closure $cl\left(L'_{n0} \right)$ of the operator $L'_{n0} $ is called the minimal operator generated by differential-operator expression \eqref{GrindEQ__2_2_} and it is denoted by $L_{n0} $ in $L_{n}^{2} ,{\rm \; }n\ge 1$.
 The operator $L_{0} $ defined by
\[D\left(L_{0} \right):=\left\{u=\left(u_{n}\right):u_{n} \in D\left(L_{n0} \right),{\rm \; }n\ge 1,{\rm \; }\sum _{n=1}^{\infty }\left\| L_{n0} u_{n} \right\| _{L_{n}^{2} }^{2} <+\infty  \right\},\]

\[L_{0} u:=\left(L_{n0}u_{n}\right),{\rm \; }u\in D\left(L_{0} \right), L_{0} :D\left(L_{0} \right)\subset L^{2} \to L^{2} \]
is called a minimal operator ( multipoint ) generated by differential-operator expression \eqref{GrindEQ__2_1_} in Hilbert space $L^{2} $ and denoted by $L_{0} =\mathop{\oplus }\limits_{n=1}^{\infty } L_{n0} $.

 In a similar way the minimal operator (twopoints) $L_{n0}^{+} $ in $L_{n}^{2} ,{\rm \; }n\ge 1$ for the formally adjoint linear differential-operator expression \eqref{GrindEQ__2_3_} can be constructed.

 In this case the operator  $L_{0}^{+} $ defined by
\[D\left(L_{0}^{+} \right):=\left\{v:=\left(v_{n} \right):v_{n} \in D\left(L_{n0}^{+} \right),{\rm \; }n\ge 1,{\rm \; }\sum _{n=1}^{\infty }\left\| L_{n0}^{+} v_{n} \right\| _{L_{n}^{2} }^{2} <+\infty  \right\},\]
$L_{0}^{+} v:=\left(L_{n0}^{+} v_{n} \right),{\rm \; }v\in D\left(L_{0}^{+} \right), L_{0}^{+} :D\left(L_{0}^{+} \right)\subset L^{2} \to L^{2}$
is called a minimal operator ( multipoint ) generated by $l^{+} \left(v\right)=(l_{n}^{+}(v_{n}))$ in the Hilbert space $L^{2} $ and denoted by $L_{0}^{+} =\mathop{\oplus }\limits_{n=1}^{\infty } L_{n0}^{+}$.

  Note that the following proposition is true.

\noindent \textbf{2.1. Theorem.} The minimal operators $L_{0} $ and $L_{0}^{+} $ are densely defined closed operators in $L^{2} $.\\

 The following defined operators in $L^{2} $

  $L:=\left(L_{0}^{+} \right)^{*} =\mathop{\oplus }\limits_{n=1}^{\infty } L_{n} $ and $L^{+} :=\left(L_{0} \right)^{*} =\mathop{\oplus }\limits_{n=1}^{\infty } L_{n}^{+} $

\noindent are called maximal operators ( multipoint ) for the differential-operator expression $l\left(\cdot \right)$ and $l^{+} \left(\cdot \right)$ respectively.It is clear that $Lu=\left(l_{n} \left(u_{n} \right)\right),{\rm \; }u\in D\left(L\right)$,

\[D\left(L\right):=\left\{u=\left(u_{n} \right)\in L^{2} :u_{n} \in D\left(L_{n} \right),{\rm \; }n\ge 1{\rm \; ,\; }\sum _{n=1}^{\infty }\left\| L_{n} u_{n} \right\| _{L_{n}^{2} }^{2} <\infty  \right\},\]

\[L^{+} v=\left(l_{n}^{+} \left(v_{n} \right)\right),{\rm \; }v\in D\left(L^{+} \right),\]

\[D\left(L^{+} \right):=\left\{v=\left(v_{n} \right)\in L^{2} :v_{n} \in D\left(L_{n}^{+} \right),{\rm \; }n\ge 1{\rm \; ,\; }\sum _{n=1}^{\infty }\left\| L_{n}^{+} v_{n} \right\| _{L_{n}^{2} }^{2} <\infty  \right\}\]
and $L_{0} \subset L$, $L_{0}^{+} \subset L^{+} $.

 Furthermore, the validity of following proposition is clear.

\noindent \textbf{2.2. Theorem.} The domain of the operators $L$ and $L_{0} $ are
\[\begin{array}{l} {D\left(L\right)=\left\{u=\left(u_{n} \right)\in L^{2} :\left(1\right){\rm \; for\; each\; }n\ge 1{\rm \; vector\; function\; }u_{n} \in L_{n}^{2} ,{\rm \; derivative\; } u'_{n} \right. } \\ {{\rm \; \; \; \; \; \; \; \; \; \; \; \; \; \; \; \; \; \; \; \; \; \; \; \; \; \; \; \; \; \; \; \; \; \; \; \; \; \; \; \; is\; absolutely\; continuous\; in\; interval\; }\Delta _{n} ;} \\ {{\rm \; \; \; \; \; \; \; \; \; \; \; \; \; \; \; \; \; \; \; \; \; \;\; \; \; \; \; \; \; \; \; \; \;\; \; \; \; \; \; }\left(2\right){\rm \; }l_{n} \left(u_{n} \right)\in L_{n}^{2} ,{\rm \; }n\ge 1;{\rm \; }\left. \left(3\right){\rm \; }l\left(u\right)=\left(l_{n} \left(u_{n} \right)\right)\in L^{2} \right\}} \\
{{\rm \; \; \; \; \; \; \; \; \; }=\left\{u=\left(u_{n} \right)\in L^{2} :u_{n} \in D\left(L_{n} \right),{\rm \; }n\ge 1{\rm \; and\; }l\left(u\right)=\left(l_{n} \left(u_{n} \right)\right)\in L^{2} \right\},} \end{array}\]

$ D(L_{0})=\left\{u=\left(u_{n} \right)\in D\left(L\right):u_{n} \left(a_{n} \right)=u_{n} \left(b_{n} \right)=u'_{n} \left(a_{n} \right)=u'_{n} \left(b_{n} \right)=0,{\rm \; }n\ge 1\right\}$.
\textbf{2.3. Remark.} If for any $n\ge 1$, $u_{n} \in D\left(L_{n0} \right)$, then $u=\left(u_{n} \right)$ may not be in $D\left(L_{0} \right)$ in general.
 Indeed, choose the function in form
\[u_{n} \left(t\right):=c_{n} sin^{2} \left(n\pi \frac{t-a_{n} }{b_{n} -a_{n} } \right)f_{n} ,{\rm \; }\]
\noindent $c_{n} \in {\mathbb C},{\rm \; }t\in \Delta _{n} ,{\rm \; }\left(f_{n} \right)\in D\left(A\right),{\rm \; \; }f_{n} \ne 0,{\rm \; }\alpha _{n} =\left\| f_{n} \right\| _{H} {\rm ,\; }n\ge 1.$

 In this case it is easily to see that $u_{n} \in D\left(L_{n0} \right),{\rm \; }n\ge 1$. On the other hand the simple calculations shows that
\[\left\| u\right\| _{L^{2} }^{2} =\sum _{n=1}^{\infty }\int _{a_{n} }^{b_{n} }\left\| u_{n} \right\| _{L_{n}^{2} }^{2} dt =\sum _{n=1}^{\infty }\alpha _{n}^{2} c_{n}^{2}   \int _{a_{n} }^{b_{n} }sin^{4} \left(n\pi \frac{t-a_{n} }{b_{n} -a_{n} } \right)dt \] \\
\[ {{\rm \; \; \; \; \; \; \; \; \; \; \; \;  \; \; \; \; \; \; \; \; \; \; \; \; \;}=\frac{3}{4} \sum _{n=1}^{\infty }\alpha _{n}^{2} c_{n}^{2} \left(b_{n} -a_{n} \right).}\]

 Here, if put $\alpha _{n} :=\frac{1}{n} ,{\rm \; }c_{n} :=\left[\frac{4}{3\left(b_{n} -a_{n} \right)} \right]^{1/2} \cdot \frac{1}{\alpha _{n} } ,{\rm \; }n\ge 1$,

\noindent then from the last relation implies that $\left\| u\right\| _{L^{2} }^{2} =\sum 1 =+\infty $, i.e $u\notin L^{2}. $

\noindent \textbf{2.4. Remark.} If $A_{n} \in B\left(H\right),{\rm \; }n\ge 1$ and $\mathop{\sup }\limits_{n\ge 1} \left\| A_{n} \right\| \le c<+\infty $, then for any $u=\left(u_{n} \right)\in L^{2} $ we have $\left(Au\right)=\left(A_{n} u_{n} \right)\in L^{2} $.

Now the following results can be proved .

\noindent \textbf{2.5. Theorem.} If\textbf{ }a minimal\textbf{ }operator $L_{0} $ is formally normal in $L^{2} $, then $D\left(L_{0} \right)\subset \mathop{W_{2}^{2} }\limits^{0} $ and  $AD\left(L_{0} \right)\subset L^{2} $.

\noindent \textbf{2.6. Theorem} If $A^{1/2} W_{2}^{2} \subset W_{2}^{2} $, then minimal operator $L_{0} $ is formally normal in $L^{2} $.

\noindent \textbf{Proof:} In this case from the following relations
\[L_{0}^{+} u=-u''-iAu=\left(-u''+iAu\right)-2iAu=L_{0} u-2iAu,{\rm \; }u\in D\left(L_{0} \right),\]
\[L_{0} u=-u''+iAu=\left(-u''-iAu\right)+2iAv=L_{0}^{+} u+2iAu,{\rm \; }u\in D\left(L_{0}^{+} \right),\]
implies that $D\left(L_{0} \right)=D\left(L_{0}^{+} \right)$. Since $D\left(L_{0}^{+} \right)\subset D\left(L_{0}^{*} \right)=D\left(L^{+} \right)$, it is obtained that $D\left(L_{0} \right)\subset D\left(L^{+} \right)$.\\

 On the other hand for any $u\in D\left(L_{0} \right)$
\[\begin{array}{l} {\left\| L_{0} u\right\| _{L^{2} }^{2} =\left(-u''+iAu,-u''+iAu\right)_{L^{2} } =\left\| u''\right\| _{L^{2} }^{2} +i\left[\left(u'',Au\right)_{L^{2} } -\left(Au,u''\right)_{L^{2} } \right]+\left\| Au\right\| _{L^{2} }^{2} } \\ {{\rm \; \; \; \; \; \; \; \; \; \; \; }=\left\| u''\right\| _{L^{2} }^{2} +\left\| Au\right\| _{L^{2} }^{2} } \end{array}\]
and
\[\begin{array}{l} {\left\| L^{+} u\right\| _{L^{2} }^{2} =\left(-u''-iAu,-u''-iAu\right)_{L^{2} } =\left\| u''\right\| _{L^{2} }^{2} -i\left[\left(u'',Au\right)_{L^{2} } -\left(Au,u''\right)_{L^{2} } \right]+\left\| Au\right\| _{L^{2} }^{2} } \\ {{\rm \; \; \; \; \; \; \; \; \; \; \; }=\left\| u''\right\| _{L^{2} }^{2} +\left\| Au\right\| _{L^{2} }^{2} .} \end{array}\]
From this it is established that operator $L_{0} $ is formally normal in $L^{2} $.

\noindent \textbf{2.7 Remark.} If  $A_{n} \in B\left(H\right),{\rm \; }n\ge 1$ and $\mathop{\sup }\limits_{n\ge 1} \left\| A_{n} \right\| \le c<+\infty $, then $D\left(L_{0} \right)=D\left(L_{0}^{+} \right)$ and $D\left(L\right)=D\left(L^{+} \right).$

\noindent \textbf{2.8 Remark.} If $AW_{2}^{2} \subset L^{2} $, then $D\left(L_{0} \right)=D\left(L_{0}^{+} \right)$ and $D\left(L\right)=D\left(L^{+} \right)$.

\section{\textbf{Description of Normal Extensions of the Minimal Operator} }

\indent In this section the main purpose is to describe all normal extensions of the minimal operator  \textbf{$L_{0} $ }in\textbf{ $L^{2} $ }in terms in the boundary values of the endpoints of the\textbf{ }subintervals .

 In first will be shown that there exists normal extension of the minimal operator \textbf{$L_{0} $}. Consider the following extension of the\textbf{ }minimal operator \textbf{$L_{0} $}
\[\left\{\begin{array}{l} {{\rm \; \; \; \; }\widetilde{L}u:=-u''+iAu,{\rm \; }AW_{2}^{2} \subset W_{2}^{2} ,} \\ {D\left(\widetilde{L}\right)=\left\{u=\left(u_{n} \right)\in W_{2}^{2} :u_{n} \left(a_{n} \right)=u_{n} \left(b_{n} \right),{\rm \; }u'_{n} \left(a_{n} \right)=u'_{n} \left(b_{n} \right),{\rm \; }n\ge 1\right\}.} \end{array}\right. \]
 Under the condition to the coefficient $A$ we have
\[\begin{array}{l} {{\left(\widetilde{L}u,v\right)_{L^{2} } =\left(-u'',v\right)_{L^{2} } +i\left(Au,v\right)_{L^{2} }}} \\
{{\rm \; \; \; \; \; \; \; }=\left(-u',v\right)_{L^{2} } ^{{'} } +\left(u,v'\right)_{L^{2} } ^{{'} } -\left(u,v''\right)_{L^{2} } +\left(u,-iAv\right)_{L^{2} } {\rm \; \; \; \; \; \; \; \; \; \; \; \; \; }} \\
{{\rm \; \; \; \; \; \; \; }=\sum _{n=1}^{\infty }\left[\left(u'_{n} \left(b_{n} \right),v_{n} \left(b_{n} \right)-v_{n} \left(a_{n} \right)\right)_{n} +\left(u_{n} \left(a_{n} \right),v'_{n} \left(b_{n} \right)-v'_{n} \left(a_{n} \right)\right)_{n} \right]} + \\
 {{\rm \; \; \; \; \; \; \; \; \; \; \; \;\; \; }+{\left(u,-v''-iAv\right)_{L^{2} } }} \end{array}\]
From this it is obtained \\
\[\left\{\begin{array}{l} {{\rm \; \; \; \; }\widetilde{L}^{*}v:=-v''-iAu,{\rm \; }AW_{2}^{2} \subset W_{2}^{2} ,} \\ {D\left(\widetilde{L}^{*}\right)=\left\{v=\left(v_{n} \right)\in W_{2}^{2} :v_{n} \left(a_{n} \right)=v_{n} \left(b_{n} \right),{\rm \; }v'_{n} \left(a_{n} \right)=v'_{n} \left(b_{n} \right),{\rm \; }n\ge 1\right\}.} \end{array}\right. \]
\noindent In this case it is clear that $D\left(\widetilde{L}\right)=D\left(\widetilde{L}^{*} \right)$. On the other hand since for each $u\in D\left(\widetilde{L}\right)$
\[\begin{array}{l} {\left\| \widetilde{L}u\right\| _{L^{2} }^{2} =\left\| u''\right\| _{L^{2} }^{2} +i\left[\left(u'',Au\right)_{L^{2} } -\left(Au,u''\right)_{L^{2} } \right]+\left\| Au\right\| _{L^{2} }^{2} ,} \\ {\left\| \widetilde{L}^{*} u\right\| _{L^{2} }^{2} =\left\| u''\right\| _{L^{2} }^{2} -i\left[\left(u'',Au\right)_{L^{2} } -\left(Au,u''\right)_{L^{2} } \right]+\left\| Au\right\| _{L^{2} }^{2} } \end{array}\]
and \\
\[\begin{array}{l} {\left(u'',Au\right)_{L^{2} } -\left(Au,u''\right)_{L^{2} } =\left(u',Au\right)_{L^{2} } ^{{'} } -\left(u,Au'\right)_{L^{2} } ^{{'} } } \\ \noindent {=\sum _{n=1}^{\infty }\left[\left(u'_{n} \left(a_{n} \right),A_{n} \left(u_{n} \left(b_{n} \right)-u_{n} \left(a_{n} \right)\right)\right)_{n} -\left(u_{n} \left(a_{n} \right),A_{n} \left(u'_{n} \left(b_{n} \right)-u'_{n} \left(a_{n} \right)\right)\right)_{n} \right] =0,} \end{array}\]
then $\left\| \widetilde{L}u\right\| _{L^{2} } =\left\| \widetilde{L}^{*} u\right\| _{L^{2} } $ for every $u\in D\left(\widetilde{L}\right)$. Consequently, $\widetilde{L}$ is a normal extension of the  minimal operator $L_{0} $.

 The following result established relationship between normal extensions of $L_{0} $ and normal extensions of  $L_{n0} ,{\rm \; }n\ge 1$.

\noindent \textbf{3.1 Theorem.} If $\widetilde{L}$ is a normal extension of the minimal operator $L_{0} $ in $L^{2} $, then for any $n\ge 1$,

\[D\left(\widetilde{L_{n} }\right)=P_{n} D\left(\widetilde{L}\right),{\rm \; }\widetilde{L_{n} }u_{n} =L_{n} u_{n} , \]
where $P_{n} :L^{2} \to L_{n}^{2} $ is an orthogonal projection, is a normal extension of the minimal operator $L_{n0} $ in $L_{n}^{2} ,{\rm \; }n\ge 1$.

\noindent \textbf{Proof:} Indeed, in this case firstly it is clear that

\textbf{ $D\left(L_{n0} \right)\subset P_{n} D\left(\widetilde{L}\right)\subset D\left(L_{n} \right),{\rm \; }n\ge 1$ }and\textbf{ $D\left(L_{n0}^{+} \right)\subset P_{n} D\left(\widetilde{L}^{*} \right)\subset D\left(L_{n}^{+} \right),{\rm \; }n\ge 1$}.  Now define
\[\widetilde{L_{n} }u_{n} :=l_{n} \left(u_{n} \right),{\rm \; }D\left(\widetilde{L_{n} }\right)=P_{n} D\left(\widetilde{L}\right),{\rm \; }n\ge 1\]
which is an extension of the minimal operator $L_{n0} $ in $L_{n}^{2} $, $n\ge 1$.

 Prove that an extension $\widetilde{L_{n} }$ is a normal in $L_{n}^{2} ,{\rm \; }n\ge 1$. First of all, note that from above writting relations and normality of $\widetilde{L}$ it is implies that $D\left(\widetilde{L_{n} }\right)=D\left(\widetilde{L_{n} }^{*} \right)$. \\
\indent Indeed, now put any $u_{n} \in D\left(\widetilde{L_{n} }\right),{\rm \; }n\ge 1$, then $u^{*} =\left\{0,0,\cdots ,u_{n} ,0,\cdots \right\}\in D\left(\widetilde{L}\right)$. Because $P_{n} u^{*} \in D\left(\widetilde{L_{n} }\right)$. On the other hand, since $D\left(\widetilde{L}\right)=D\left(\widetilde{L}^{*} \right)$ , then $u^{*} \in D\left(\widetilde{L}^{*} \right)$. Hence, $u_{n} \in D\left(\widetilde{L_{n} }^{*} \right),{\rm \; }n\ge 1$. From this, it is clear that $D\left(\widetilde{L_{n} }\right)\subset D\left(\widetilde{L_{n} }^{*} \right)$.

 A similar way it is shown that $D\left(\widetilde{L_{n} }^{*} \right)\subset D\left(\widetilde{L_{n} }\right)$. Finally we have $D\left(\widetilde{L_{n} }\right)=D\left(\widetilde{L_{n} }^{*} \right)$.

 On the other hand for any $u_{n} \in D\left(\widetilde{L_{n} }\right),{\rm \; }n\ge 1$ the function $u^{*} \in D\left(\widetilde{L}\right)$. Since $\left\| \widetilde{L}u^{*} \right\| _{L^{2} } =\left\| \widetilde{L}^{*} u^{*} \right\| _{L^{2} } $, then \\
\[\left\| \widetilde{L_{n} }u_{n} \right\| _{L_{n}^{2} } =\left\| \widetilde{L_{n} }^{*} u_{n} \right\| _{L_{n}^{2} } ,{\rm \; }u_{n} \in D\left(\widetilde{L_{n} }\right),{\rm \; }n\ge 1.\]
\indent Consequently, it is clear that extension $\widetilde{L_{n} }$ is a normal in $L_{n}^{2} $, $n\ge 1$.

 The minimal operator $L_{n0}^{r} $ generated by differential expression $l_{n}^{r} \left(u_{n} \right):=-u''_{n} \left(t\right)$ in Hilbert space $L_{n}^{2} ,{\rm \; }n\ge 1$ is symmetric and has equal defect indexes (dim$H$, dim$H$). Then the minimal operator $L_{n0}^{r} $ in $L_{n}^{2} ,{\rm \; }n\ge 1$ has a space of boundary values $\left({\rm {\mathfrak H}}_{n} ,\gamma _{1}^{\left(n\right)} ,\gamma _{2}^{\left(n\right)} \right),{\rm \; }n\ge 1$ \cite{Gor} and one of these spaces is in following form
\[{\rm {\mathfrak H}}_{n} =H_{n} \oplus H_{n} ,{\rm \; }\gamma _{1}^{\left(n\right)} \left(u_{n} \right)=\left\{-u_{n} \left(a_{n} \right),{\rm \; }u_{n} \left(b_{n} \right)\right\},{\rm \; }\gamma _{2}^{\left(n\right)} \left(u_{n} \right)=\left\{u'_{n} \left(a_{n} \right),{\rm \; }u'_{n} \left(b_{n} \right)\right\},\]
$u_{n} \in D\left(L_{n}^{r} \right)$, where $L_{n}^{r} $ is defined a maximal operator generated by differential expression $l_{n}^{r} \left(\cdot \right)=-d^{2} /dt^{2} $ in the space $L_{n}^{2} ,{\rm \; }n\ge 1$.

 Now we can prove the following main result of this section in which is given a description of all normal extension  of the minimal operator $L_{0} $ in $L^{2} $ in terms of boundary values of vector functions at the endpoints of subintervals.

\noindent \textbf{3.2 Theorem.} Let $A^{1/2} W_{2}^{2} \subset W_{2}^{2} $. If $\widetilde{L}=\mathop{\oplus }\limits_{n=1}^{\infty } \widetilde{L_{n} }$ is a normal extension of the minimal operator $L_{0} $ in $L^{2} $, then it is generated by differential-operator expression $l\left(\cdot \right)=\left(l_{n} \left(\cdot \right)\right)$ with boundary conditions

\[\left(W_{n} -E\right)\gamma _{1}^{\left(n\right)} \left(u_{n} \right)+i\left(W_{n} +E\right)\gamma _{2}^{\left(n\right)} \left(u_{n} \right)=0,{\rm \; }u_{n} \in D\left(L_{n} \right),\]
where, $\left({\rm {\mathfrak H}}_{n} ,\gamma _{1}^{\left(n\right)} ,\gamma _{2}^{\left(n\right)} \right)$ is a space of boundary values of the minimal operator $L_{n0}^{r} $ in $L_{n}^{2} $,$W_{n} $ and $\left(\begin{array}{cc} {A_{n}^{1/2} } & {0} \\ {0} & {A_{n}^{1/2} } \end{array}\right)W_{n} \left(\begin{array}{cc} {A_{n}^{-1/2} } & {0} \\ {0} & {A_{n}^{-1/2} } \end{array}\right)$ are unitary operators in ${\rm {\mathfrak H}}_{n} ,{\rm \; }n\ge 1$. The unitary operator $W=\mathop{\oplus }\limits_{n=1}^{\infty } W_{n} $ in ${\rm {\mathfrak H}}=\mathop{\oplus }\limits_{n=1}^{\infty } {\rm {\mathfrak H}}_{n} $ is determined uniquely by the extension $\widetilde{L}$, i.e $\widetilde{L}=L_{W} $.

\noindent \textbf{Proof:} In this case by the Theorem 3.1 $\widetilde{L_{n} }$ is a normal extension of the minimal operator $L_{n0} $ in $L_{n}^{2} ,{\rm \; }n\ge 1$. On the other hand, it is clear that $cl\left(Re\widetilde{L_{n} }\right)$ is a some selfadjoint extension of the minimal operator $L_{n0}^{r} $ generated by the differential expression $l_{n}^{r} \left(\cdot \right)=-d^{2} /dt^{2} $ in the $L_{n}^{2} ,{\rm \; }n\ge 1$. In this case there exists space of boundary values $\left({\rm {\mathfrak H}}_{n} ,\gamma _{1}^{\left(n\right)} ,\gamma _{2}^{\left(n\right)} \right)$ for the $L_{n0}^{r} $ in $L_{n}^{2} ,{\rm \; }n\ge 1$ \cite{Gor}. Therefore, the selfadjoint extension $cl\left(Re\widetilde{L_{n} }\right)$ in $L_{n}^{2} $ is generated by the differential expression $l_{n}^{r} \left(\cdot \right)=-d^{2} /dt^{2} $ and boundary condition

\begin{equation} \label{GrindEQ__3_1_}
\left(W_{n} -E\right)\gamma _{1}^{\left(n\right)} \left(u_{n} \right)+i\left(W_{n} +E\right)\gamma _{2}^{\left(n\right)} \left(u_{n} \right)=0,{\rm \; }u_{n} \in D\left(L_{n} \right),
\end{equation}
where $W_{n} $ is a unitary operator in ${\rm {\mathfrak H}}_{n} ,{\rm \; }n\ge 1$. Note that the unitary operator $W_{n} $ is determined uniquely by the extension $cl\left(Re\widetilde{L_{n} }\right)$, i.e. $cl\left(Re\widetilde{L_{n} }\right)=L_{n} \left(W_{n} \right),{\rm \; }n\ge 1$ \cite{Gor}, \cite{RoKh}.

 On the other hand, $cl\left(Im\widetilde{L_{n} }\right)$ is a selfadjoint operator which is acting in $L_{n}^{2} $ with domain $D\left(cl\left(Re\widetilde{L_{n} }\right)\right),{\rm \; }n\ge 1$. Since selfadjoint operators $cl\left(Re\widetilde{L_{n} }\right)$ and $cl\left(Im\widetilde{L_{n} }\right)$ are commutative in the space $L_{n}^{2} ,{\rm \; }n\ge 1$, then for \\
\noindent each  $u_{n} \in D\left(cl\left(Re\widetilde{L_{n} }\right)\right)$

\noindent ${\left(cl\left(Re\widetilde{L_{n} }\right)u_{n} ,cl\left(Im\widetilde{L_{n} }\right)u_{n} \right)_{L_{n}^{2} } -\left(cl\left(Im\widetilde{L_{n} }\right)u_{n} ,cl\left(Re\widetilde{L_{n} }\right)u_{n} \right)_{L_{n}^{2} }}=$ \\
\noindent $=\left(-u''_{n} ,A_{n} u_{n} \right)_{L_{n}^{2} } -\left(A_{n} u_{n} ,-u''_{n} \right)_{L_{n}^{2} } $ \\
\noindent $=\left(\gamma _{1}^{\left(n\right)} \left(A_{n}^{1/2} u_{n} \right),\gamma _{2}^{\left(n\right)} \left(A_{n}^{1/2} u_{n} \right)\right)_{{\rm {\mathfrak H}}_{n} } -\left(\gamma _{2}^{\left(n\right)} \left(A_{n}^{1/2} u_{n} \right),\gamma _{1}^{\left(n\right)} \left(A_{n}^{1/2} u_{n} \right)\right)_{{\rm {\mathfrak H}}_{n} } =0$

Since for each $n\ge 1$ $A_{n}^{-1} \in B\left(H\right)$, then last equation means that the linear relation
\[\widetilde{\theta _{n} }:=\left\{\left\{\gamma _{1}^{\left(n\right)} \left(A_{n}^{1/2} u_{n} \right),\gamma _{2}^{\left(n\right)} \left(A_{n}^{1/2} u_{n} \right)\right\}:u_{n} \in D\left(cl\left(Re\widetilde{L_{n} }\right)\right)\right\}\]
is selfadjoint in ${\rm {\mathfrak H}}_{n} \oplus {\rm {\mathfrak H}}_{n} $ (see \cite{Cod}, \cite{RoKh}). Consequently, there is a uniquely $V_{n} $ unitary operator in ${\rm {\mathfrak H}}_{n} $, such that
\begin{equation} \label{GrindEQ__3_2_}
\left(V_{n} -E\right)\gamma _{1}^{\left(n\right)} \left(A_{n}^{1/2} u_{n} \right)+i\left(V_{n} +E\right)\gamma _{2}^{\left(n\right)} \left(A_{n}^{1/2} u_{n} \right)=0
\end{equation}\\
$ u_{n} \in D\left(cl\left(Re\widetilde{L_{n} }\right)\right),{\rm \; }n\ge 1 \cite{RoKh}$.

\noindent Now consider the mappings $\gamma _{1}^{\left(n\right)} \left(A_{n}^{1/2} \cdot \right)$ and $\gamma _{2}^{\left(n\right)} \left(A_{n}^{1/2} \cdot \right)$ defined in ${\rm {\mathfrak H}}_{n} ,{\rm \; }n\ge 1$. For $u_{n} \in D\left(cl\left(Re\widetilde{L_{n} }\right)\right),{\rm \; }n\ge 1$, it is clear that

\[\gamma _{1}^{\left(n\right)} \left(A_{n}^{1/2} u_{n} \right)=\left\{-A_{n}^{1/2} u_{n} \left(0\right),A_{n}^{1/2} u_{n} \left(1\right)\right\}\]

\[=\left(\begin{array}{cc} {A_{n}^{1/2} } & {0} \\ {0} & {A_{n}^{1/2} } \end{array}\right)\left(\begin{array}{c} {-u_{n} \left(0\right)} \\ {u_{n} \left(1\right)} \end{array}\right)=\left(\begin{array}{cc} {A_{n}^{1/2} } & {0} \\ {0} & {A_{n}^{1/2} } \end{array}\right)\gamma _{1}^{\left(n\right)} \left(u_{n} \right)\]
and

\[\begin{array}{l} {\gamma _{2}^{\left(n\right)} \left(A_{n}^{1/2} u_{n} \right)=\left\{A_{n}^{1/2} u'_{n} \left(0\right),A_{n}^{1/2} u'_{n} \left(1\right)\right\}} \\ {{\rm \; \; \; \; \; \; \; \; \; \; \; \; \; \; \; \; \; \; }=\left(\begin{array}{cc} {A_{n}^{1/2} } & {0} \\ {0} & {A_{n}^{1/2} } \end{array}\right)\left(\begin{array}{c} {u'_{n} \left(0\right)} \\ {u'_{n} \left(1\right)} \end{array}\right)=\left(\begin{array}{cc} {A_{n}^{1/2} } & {0} \\ {0} & {A_{n}^{1/2} } \end{array}\right)\gamma _{2}^{\left(n\right)} \left(u_{n} \right).} \end{array}\]
From the last result and \eqref{GrindEQ__3_2_} it is obtained that

\[\left(V_{n} -E\right)\left(\begin{array}{cc} {A_{n}^{1/2} } & {0} \\ {0} & {A_{n}^{1/2} } \end{array}\right)\gamma _{1}^{\left(n\right)} \left(u_{n} \right)+i\left(V_{n} +E\right)\left(\begin{array}{cc} {A_{n}^{1/2} } & {0} \\ {0} & {A_{n}^{1/2} } \end{array}\right)\gamma _{2}^{\left(n\right)} \left(u_{n} \right)=0,{\rm \; }n\ge 1.\]

 On the other hand, since $A_{n} >0,{\rm \; }n\ge 1$, then the last relation is equivalent to the following equation \\
\noindent $\left(\left(\begin{array}{cc} {A_{n}^{-1/2} } & {0} \\ {0} & {A_{n}^{-1/2} } \end{array}\right)V_{n} \left(\begin{array}{cc} {A_{n}^{1/2} } & {0} \\ {0} & {A_{n}^{1/2} } \end{array}\right)-E\right)\gamma _{1}^{\left(n\right)} \left(u_{n} \right)$ +\\
{\rm \; \;\;\;}$+i\left(\left(\begin{array}{cc} {A_{n}^{-1/2} } & {0} \\ {0} & {A_{n}^{-1/2} } \end{array}\right)V_{n} \left(\begin{array}{cc} {A_{n}^{1/2} } & {0} \\ {0} & {A_{n}^{1/2} } \end{array}\right)-E\right)\gamma _{2}^{\left(n\right)} \left(u_{n} \right)=0,$ \\
\noindent $u_{n} \in D\left(cl\left(Re\widetilde{L_{n} }\right)\right),{\rm \; }n\ge 1.$ \\
 From the last relation, \eqref{GrindEQ__3_1_} and uniqueness of operator $W_{n} $ it is clear that
\[W_{n} =\left(\begin{array}{cc} {A_{n}^{-1/2} } & {0} \\ {0} & {A_{n}^{-1/2} } \end{array}\right)V_{n} \left(\begin{array}{cc} {A_{n}^{1/2} } & {0} \\ {0} & {A_{n}^{1/2} } \end{array}\right), \]
i.e operator

\[V_{n} =\left(\begin{array}{cc} {A_{n}^{1/2} } & {0} \\ {0} & {A_{n}^{1/2} } \end{array}\right)W_{n} \left(\begin{array}{cc} {A_{n}^{-1/2} } & {0} \\ {0} & {A_{n}^{-1/2} } \end{array}\right):{\rm {\mathfrak H}}_{n} \to {\rm {\mathfrak H}}_{n} ,{\rm \; }n\ge 1\]
must be unitary operator.

Finally, we give one result on the spectrum structure of the normal extensions.\\

\noindent \textbf{3.3 Theorem.} For the spectrum of the normal extension $\widetilde{L}=\mathop{\oplus }\limits_{n=1}^{\infty } \widetilde{L_{n} }$ of the minimal operator $L_{0} $ in the space $L^{2} $ the following formulas are true

\[\sigma _{p} \left(\widetilde{L}\right)=\bigcup _{n=1}^{\infty }\sigma _{p} \left(\widetilde{L_{n} }\right) ,{\rm \; \; }\bigcap _{n=1}^{\infty }\sigma _{c} \left(A_{n} \right)\subset  {\rm \; }\sigma _{c} \left(A\right)\subset \bigcup _{n=1}^{\infty }\sigma _{c} \left(A_{n} \right) .\]
\textbf{3.4 Remark.} In this work for the simplicity of explanation the multipoint differential-operator expression has been considered in form \eqref{GrindEQ__2_1_}-\eqref{GrindEQ__2_2_}. However using the established results in \cite{Ism2} the analogous claims can be obtained in the case when operator coefficient in \eqref{GrindEQ__2_1_}-\eqref{GrindEQ__2_2_} are any normal operator in $H_{n} ,{\rm \; }n\ge 1$.

\noindent \textbf{}
\bibliographystyle{amsalpha}

\end{document}